\documentclass[10pt]{article}
\usepackage{theorem}
\usepackage{amssymb}

\usepackage{amsmath}   
\usepackage{amsfonts}

\def\RR{\mathbb R}

\def\Z{\mathbb Z}
\def\EE{\mathbb E}
\def\FF{\mathcal F}

\newtheorem{theorem}{Theorem}[section]
\newtheorem{lemma}[theorem]{Lemma}
\newtheorem{proposition}[theorem]{Proposition}
\newtheorem{corollary}[theorem]{Corollary}
 {\theorembodyfont{\rmfamily}

\begin{document}

\title{ $5n$ Minkowski symmetrizations suffice to arrive at an
approximate Euclidean ball } 

\author{ B. Klartag\thanks{ Supported by the Israel Science Foundation
founded by the Academy of Sciences and Humanities. } \\
School of Mathematical Sciences, \\
Tel Aviv University, \\
Tel Aviv 69978, Israel }

\date{}

\maketitle              

\begin{abstract}
 This paper proves that for every convex body in $\RR^n$ there exist
 $5n-4$ Minkowski symmetrizations, which transform the body into an
 approximate Euclidean ball. This result complements the sharp $c n
 \log n$ upper estimate by J. Bourgain, J. Lindenstrauss and
 V.D. Milman, of the number of random Minkowski symmetrizations
 sufficient for approaching an approximate Euclidean ball. 
\end{abstract}

\section{Introduction}

Let $K$ be a compact convex set in $\RR^n$ and let $u$ be any
vector in $S^{n-1} = \{u ; |u| = 1\} $ where $| \cdot |$ denotes
the standard Euclidean norm in $\RR^n$. Denote by $\pi_{u} \in
O(n)$ the reflection with respect to the hyperplane through the
origin orthogonal to $u$, i.e. $ \pi_{u} x = x - 2 \langle x,u
\rangle u$.

\medskip Minkowski symmetrization (often referred to as Blaschke
symmetrization) of $K$ with respect to $u$ is defined to be the
convex set $\frac{1}{2}(\pi_{u}K + K)$. Denote by $\| \cdot \|^*$
the dual norm to $K$ (i.e. $\| x \|^* = \sup_{y \in K} \langle x,
y \rangle$). Despite the fact that $K$ is not necessarily
centrally symmetric and $\| \cdot \|^*$ need not be a norm, this
convenient notation will be used for readability. Denote by
$M^*(K)$ the half mean width of $K$, defined as $M^{*}(K) :=
\int_{S^{n-1}} \| x \|^* d \sigma (x) $, where $ \sigma $ is the
normalized rotation invariant measure on $S^{n-1}$, and $\| \cdot
\|^*$ is the dual norm.

\medskip It is easily verified that $M^*(K) =
M^*(\frac{1}{2}(\pi_{u}K + K))$, so the mean width is preserved
under Minkowski symmetrizations. Since successive Minkowski
symmetrizations make the body more symmetric in some sense, one
might expect convergence to a ball of radius $M^*(K)$.

Surprisingly, very few symmetrizations are sufficient for this
convergence; In \cite{BLM} it is proven that $c n \log n$ random
symmetrizations suffice to obtain from any convex body, a new
body $\tilde{K}$, such that $\frac{1}{2} M^* D \subset \tilde{K}
\subset 2 M^*D$ with high probability, where $D = \{u ; |u| \leq
1\}$ is the standard Euclidean ball in $\RR^n$.

\medskip The proof in \cite{BLM} can be slightly refined, and
rather than an estimate of $c n \log n$ for all bodies, in fact $c
n \log{ \frac{2 diam(K)}{M^*(K)}}$ symmetrizations are enough.
This quantity is always smaller then $c n \log n$ but in some
cases there is a substantial improvement; For example, the
$n$-dimensional cube needs only $c n$ random symmetrizations to be
transformed into an almost Euclidean ball.

\medskip In \cite{K} it was proven that the aforementioned estimate is
very tight and is actually a formula, as follows: For every convex
body $K$ at least $\tilde{c} n \log{ \frac{diam(K)}{2 M^*(K)}}$ random
symmetrizations are necessary in order for the body to become
close to a Euclidean ball. Hence, bodies such as $B(l_1^n)$ - the
$n$ dimensional cross polytope - in fact require at least $c n
\log n$ random symmetrizations.

\medskip This paper shows that there exist symmetrizations which are
better than random ones. There is a specific choice of $5n-4$
symmetrizations that transform any convex body into an approximate
Euclidean ball. The basic idea underlying the construction is
changing the notion of randomness; Rather than symmetrizing with
respect to random vectors, symmetrizations with respect to the vectors
of a random orthogonal basis will be performed at each iteration.

Six iterations of this kind suffice (totaling $6n-5$
symmetrizations \footnotemark[1]), however the role of each 
iteration is slightly different. Precisely, for the first
iteration any orthogonal basis is adequate. The remaining five
iterations are required to be with respect to random independent
orthogonal bases, and the results hold with large probability that
tends to 1 when the dimension $n$ approaches infinity.

\footnotetext[1]{ It seems at first, that six iterations consist of $6n$
symmetrizations; However, after the first iteration, the body
becomes centrally symmetric. Following that stage, the last vector
in each orthogonal basis is unnecessary, because symmetrizing
with respect to that vector would not affect the body. }

\medskip There exists a very similar symmetrization process
that leads to a slightly better estimate, and consists of $5n-4$
symmetrizations (only $4n-4$ symmetrizations, if the body is
already unconditional). This process uses symmetrizations with
respect to five orthogonal bases, some of which need not be random. An
additional basis will be used in this process, and 
will be referred to here as a Walsh basis. It actually coincides
with the regular Walsh basis for dimensions which are powers of
two. Let us describe the $5n-4$ symmetrizations process: The first
basis is chosen to be any orthogonal basis, and is used only to create
unconditionality. The second basis can be a random basis or a Walsh
basis (with respect to the first), and the corresponding
symmetrization reduces the diameter of the body
to a level of $\log n$ times its mean width. The third
basis is a Walsh basis with respect to the previous, and reduces the
diameter further, to a level of $\log \log n$ times the mean
width. The fourth basis must be, in this proof, a random orthogonal
basis and the fifth, either a Walsh basis with respect to the fourth,
or a random basis. Once the diameter is small enough, the last two
bases together transform the body to an approximate Euclidean ball.

\medskip The proof outlined below is mainly concerned with the
first process described (which is purely random). Results for the
second process are analogous to those of the first, and may be
concluded based on remarks throughout the proof.

\medskip The symbols $c, C, c^{\prime}, \tilde{c}$ denote numerical
constants which are not necessarily identical throughout this text.

\section{First Step: Initial Symmetrizations}

Let $K$ be an arbitrary convex body in $\RR^n$. For the purpose of
normalization, assume $M^*(K) = 1$. Take any orthogonal basis
$\{e_1,..,e_n\}$ and symmetrize $K$ with respect to the vectors
$e_1, .., e_n$ to obtain the new body $\tilde{K}$. Since
orthogonal reflections commute, $\tilde{K}$ is invariant under
reflection with respect to $e_i$, for $1 \leq i \leq n$. Therefore
$\tilde{K}$ is unconditional with respect to the basis
$\{e_1,..,e_n\}$. By Lemma 3.2 from \cite{K} there exists a
universal constant $c$ such that, 
$$ \tilde{K} \ \subset \ c \sqrt{n} \ conv \{\pm e_i\}_{i=1}^n = c \sqrt{n}
B(l_1^n) $$ 
A specific body will be referred to in this section: $Q = \sqrt{n} \ 
conv \{\pm e_i\}_{i=1}^n$. After a certain symmetrization process its
diameter decays from $\sqrt n$ to $\tilde{c} \log n$ with high
probability. Clearly, applying the same set of symmetrizations to
$\tilde{K}$ will reduce its diameter to less than $\tilde{c} \log n$. 
 
\begin{proposition}
Let $\{e_1,..,e_n\}$ be an orthogonal basis in $\RR^n$, and let $Q
= \sqrt{n} \ conv \{\pm e_i\}_{i=1}^n$. Let $\mu_n$ be the unique
rotation invariant probability measure on $O(n)$. Suppose that $\{
u_1, .., u_n \} \in O(n)$ is chosen randomly, according to $\mu_n$. After
symmetrizing $Q$ with respect to $u_1,..,u_{n-1}$ a new body
$\tilde{Q}$ is obtained. \\
Claim:
$$ diam(\tilde{Q}) \leq c \log n $$
with probability greater than $1 - \frac{1}{n^{10}}$.
\label{cross_polytope}
\end{proposition}

Remark: The number `10' in the expression $1 - \frac{1}{n^{10}}$
is of course arbitrary, and may be replaced by any other constant.
Such a replacement will influence the constant `c' in the
concluded inequality ``$ diam(\tilde{Q}) \leq c \log n $''.

\begin{corollary}
For every convex body $K \subset \RR^n$ with $M^*(K) = 1$, there
exist $2n-1$ symmetrizations which transform $K$ into $\tilde{K}$,
where $diam(\tilde{K}) < c \log n$ and $\tilde{K}$ is
unconditional with respect to some orthogonal basis.
\label{corollary_1}
\end{corollary}

Following is a simple and well-known lemma. For completeness it
will be proven at the end of this section.

\begin{lemma}
Let $\{e_i\}_{i=1}^n$ be any orthogonal basis, and let
$\{u_i\}_{i=1}^n$ be a random orthogonal basis. Then for all $1
\leq  i,j \leq n$:
$$ |\langle u_i, e_j \rangle| \leq c_1 \frac{\sqrt{\log
n}}{\sqrt n}$$
with probability greater than $1 - \frac{1}{n^{10}}$.
\label{maximum_of_psi_2}
\end{lemma}

\bigskip \emph{Proof of Proposition \ref{cross_polytope}}: Denote by
$\| \cdot \|$ the dual norm of $Q$ (i.e. $\| x \| = \sup_{y \in Q}
\langle x,y \rangle = \sqrt n \max_{i} |\langle x, e_i \rangle|
$). The dual norm of $\tilde{Q}$ is, by definition (recall that the
Minkowski sum of bodies is equivalent to the sum of their dual
norms):
$$ \big \| \big | x \big \| \big | = \frac{1}{2^{n-1}} \sum_{D \subset
\{1,..,n-1\} } \big\| (\prod_{i \in D} \pi_{u_i}) x \big \| =
\frac{1}{2^{n-1}} \sum_D \sqrt n \max_j \big| \langle \prod_{i \in D}
\pi_{u_i} x, e_j \rangle\big| $$
Substitute $x = \sum_i \langle x, u_i \rangle u_i $. Since
reflecting with respect to $u_i$ means switching the sign of the
$i^{th}$ coordinate in $\{u_1, .., u_n\}$ basis,
$$ \big \| \big | x \big \| \big | = \EE_{\varepsilon}
\max_j \big | \sqrt n \sum_i \varepsilon_i \langle x, u_i \rangle
\langle u_i, e_j \rangle | $$
where $\varepsilon = (\varepsilon_i)_{i=1}^n$ is uniformly distributed
in $\{ \pm 1\}^n$. Therefore, $\|| x \||$ is the expectation of a
maximum of $n$ random variables. Denote: $$ f_x^j(\varepsilon) =
\sqrt n | \sum_i \varepsilon_i \langle x, u_i \rangle \langle u_i,
e_j \rangle | $$ Then $ \|| x \|| = \EE_{\varepsilon} [ \max_j f^j_x(
\varepsilon ) ] $. For $1 \leq \alpha \leq 2$, and for any measurable
$f:\Omega \rightarrow \RR$ define $ \| f \|_{\psi_{\alpha}} = \inf \{ 
\lambda > 0 : \int_\Omega e^{|\frac{f}{\lambda}|^{\alpha}} \leq 2 \} $. 
The equivalent definitions are frequently used:
$$\|f\|_{\psi_{\alpha}} < c \ \Leftrightarrow \ (\EE|f|^p)^{\frac{1}{p}} <
c^{\prime} p^{\frac{1}{\alpha}} \ \Leftrightarrow \ Prob \{ |f| > t \} <
e^{-c^{\prime \prime} t^{\alpha}} $$
Khinchine inequality shows that the $\psi_2$ norm of $f_x^j$ is
bounded, as follows: 
$$\|f^j_x\|_p = \Big(\EE_\varepsilon \big|\sqrt n \sum_i \varepsilon_i
\langle x, u_i \rangle \langle u_i, e_j \rangle \big|^p
\Big)^{\frac{1}{p}} \leq c \sqrt p \sqrt n \sqrt{\sum_i ( \langle
x, u_i \rangle \langle u_i, e_j \rangle )^2} $$
By Lemma \ref{maximum_of_psi_2} with large probability, $|\langle
u_i, e_j \rangle| \leq c_1 \frac{\sqrt{\log n}}{\sqrt n} $. Hence,
with high probability,
$$ \|f^j_x\|_p \leq c \sqrt p \sqrt{\log n} |x| \ \ \
\Rightarrow \ \ \ \|f^j_x\|_{\psi_2} \leq c^{\prime} \sqrt{\log n}|x| $$
Since $\||x\|| = \EE_{\varepsilon} [ \max_j f^j_x( \varepsilon )
]$, the well-known estimate for the expectation of a
maximum of $\psi_2$ variables can be used (e.g. \cite{Tal} page 79, or
the remark after lemma \ref{prob_lemma} in this paper):
$$ \forall x \in \RR^n \ \||x\|| \leq c \sqrt{\log n} (c^{\prime}
\sqrt{\log n} |x|) = c \log n |x| $$ 
Thus the proposition is proven. \hfill $\square$

\medskip \emph{Remark}: For every dimension, there exists an orthogonal
basis $\{u_i\}_{i=1}^n$ such that $\forall i,j$:
$$ |\langle u_i, e_j \rangle| \leq \frac{2}{\sqrt n} $$
Such a basis is called in this paper a ``Walsh'' basis. Indeed, for
dimension $n=2^k$, the regular Walsh basis is satisfactory, 
while for other dimensions, an appropriate basis may be constructed
using sines and cosines (this basis consists of orthogonal vectors
resembling the complex valued characters of the group $\Z / n
\Z$). Instead of using Lemma \ref{maximum_of_psi_2} in the proof of
Proposition \ref{cross_polytope}, one can replace the 
random basis with a Walsh basis, obtaining yet a slightly better
result, with ``$\log n$'' replaced by ``$\sqrt{\log n}$'' in the
conclusion of Proposition \ref{cross_polytope}.

\bigskip \emph{Proof of Lemma \ref{maximum_of_psi_2}:} Since for
every $i$ the vector $u_i$ distributes uniformly over the sphere,
by the standard concentration inequality on the sphere (e.g.
first pages of \cite{MS}):
\begin{equation}
Prob \{ |\langle u_i, e_j \rangle| > \varepsilon \} \leq \sqrt
\frac{\pi}{2} e^{-\frac{\varepsilon^2 n}{2}}
\label{cap_size}
\end{equation}
Select $c_1$ (i.e. $c_1=5$) such that for $\varepsilon = c_1
\frac{\sqrt {\log n}}{\sqrt n}$, the probability in
(\ref{cap_size}) is less than $\frac{1}{n^{12}}$. Therefore, the
probability that $|\langle u_i, e_j \rangle| < c_1
\frac{\sqrt{\log n}}{\sqrt n}$ holds for all $1 \leq i,j \leq n$
is greater than $1 - \frac{1}{n^{10}}$. \hfill $\square$

\section{Second Step: Logarithmic Decay of the Diameter}

In the second step, symmetrizations will be performed with respect
to two random orthogonal bases; This section proves that this step
reduces the diameter of the body logarithmically: from $c \log n$
to $C \log {\log n}$, with probability close to $1$. Therefore, after
the second step (and a total of $4 n$ symmetrizations) the diameter is
less than $C \log {\log n}$. This proof extends that of the former section.

\medskip Let $K$ be the convex body obtained from the first step of
symmetrizations. According to Corollary \ref{corollary_1}, $M^*(K)
= 1$, $diam(K) < c \log n$, and $K$ is unconditional with respect
to some orthogonal basis (re-denote this basis as
$\{e_1,..,e_n\}$). Once again, by Lemma 3.2 from \cite{K},
$$ K \ \subset \ c \sqrt{n} \ conv \{\pm e_i\}_{i=1}^n \ = \ c
\sqrt{n} B(l_1^n) $$ 
Set $t = diam(K) < c \log n$. Clearly, $K \subset \sqrt n
B(l_1^n) \bigcap t B(l_2^n)$. As in the first step, rather
than working directly with the body $K$, symmetrize $K_t =
\sqrt n B(l_1^n) \bigcap t B(l_2^n)$. 

\begin{proposition}
Let $K_t = \sqrt{n} B(l_1^n) \bigcap t B(l_2^n)$. Assume
that $\{u_1,..,u_n\} \in O(n)$ and $\{v_1,..,v_n\} \in O(n)$ are
chosen uniformly and independently. After symmetrizations with
respect to $u_1,..,u_{n-1}$ and $v_1,..,v_{n-1}$, a new body
$\tilde{K_t}$ is obtained such that:
$$ \tilde{K_t} \ \subset \ C \log t B(l_2^n) $$
with probability greater than $1 - e^{-c\sqrt{n}}$ of choosing the
orthogonal bases.
\label{log_decay}
\end{proposition}

\begin{corollary}
For every convex body $K \subset \RR^n$ with $M^*(K) = 1$, there
exist $4n-3$ symmetrizations which transform $K$ into $\tilde{K}$,
where $diam(\tilde{K}) < c \log{\log n}$.
\label{corollary_2}
\end{corollary}

Begin by describing the body $K_t$ through its dual norm.
Denote by $\| \cdot \|^{\prime}_t$ the norm:
$$ \|x\|^{\prime}_t = \inf \{ \|x^{\prime}\|_2 + t\|x^{\prime
  \prime}\|_\infty : x = x^{\prime} + x^{\prime \prime} \} $$
The dual norm of $K_t$ is exactly $t \| \cdot
\|^{\prime}_{\frac{\sqrt{n}}{t}}$, as can be verified. Put
$(a_i^*)_{i=1}^n$ for the non-increasing rearrangement of the absolute
values of $(a_i)_{i=1}^n$. The following two lemmas are
well-known. The first lemma essentially appears in \cite{BL},
but for lack of concise references, attached here are the short
elementary proofs.

\begin{lemma}
$$ \forall x \in \RR^n \ \ \ \|x\|^{\prime}_k \approx \sqrt{\sum_{i=1}^{k^2}
(x_i^*)^2} $$ 
(and the equivalence constant is not more than $\sqrt 2$). 
\label{convolution_formula}
\end{lemma}

\emph{Proof:} For $i$ where $|x_i| \geq x_{k^2}^*$ set $x^{\prime}_i = (x_i
- sgn(x_i) x_{k^2}^*)$. For other $i$'s, set $x^{\prime}_i = 0$. Let
$x^{\prime \prime} = x - x^{\prime}$. Then:
\begin{eqnarray*}
 \lefteqn{ \|x \|^{\prime}_k \leq \|x^{\prime}\|_2 + k \|x^{\prime
 \prime}\|_\infty} \\ 
 & & = \sqrt{\sum_{i=1}^{k^2} (x_i^* - x_{k^2}^*)^2} + k x_{k^2}^* \\
 & & \leq \sqrt{2} \sqrt{\sum_{i=1}^{k^2} \big [ (x_i^* - x_{k^2}^*)^2 +
 (x_{k^2}^*)^2 \big ] } \\ 
 & & \leq \sqrt{2} \sqrt{\sum_{i=1}^{k^2} (x_i^*)^2}
\end{eqnarray*}
On the other hand, assume $x = x^{\prime} + x^{\prime \prime}$. Surely
$ x_i^* \leq x^{\prime *}_i + x^{\prime \prime *}_1 $, so:
\begin{eqnarray*}
\lefteqn{ \sqrt{\sum_{i=1}^{k^2} (x_i^*)^2} \leq \sqrt{\sum_{i=1}^{k^2} (
  x^{\prime *}_i )^2 }+ \sqrt{\sum_{i=1}^{k^2} (x^{\prime \prime *}_1)^2 }
  } \\ 
& & \leq \|x^{\prime}\|_2 + k \|x^{\prime \prime}\|_\infty
\end{eqnarray*}

\begin{lemma}
\label{prob_lemma}
Let $(X_i)_{i=1}^n$ be $\psi_1$ random variables (i.e. random
variables that satisfy: $\EE e^{|X_i|} \leq C$), and let
$(X_i^*)_{i=1}^n$ be the non-increasing rearrangement of the
$X_i$'s. Then:
$$ \EE \sqrt{ \frac{1}{k} \sum_{i=1}^k (X_i^*)^2 }\leq c_2 \log
{ \frac{2n}{k} } $$ 
\end{lemma}

\emph{Proof:} Since the $X_i$'s are $\psi_1$ variables,
\begin{equation}
\EE \frac{1}{k} \sum_{i=1}^k e^{X_i^*} \leq \EE \frac{1}{k}
\sum_{i=1}^n e^{|X_i|} \leq C \frac{n}{k} 
\label{k_max_eq}
\end{equation}
Let $(a_i)_{i=1}^k$ be any real numbers such that
$\forall i \ a_i \geq 1$. Since the function $e^{\sqrt{x}}$ is convex on $[1,
\infty)$, by Jensen inequality:
\begin{equation}
e^{\sqrt{ \frac{1}{k} \sum_{i=1}^k (a_i)^2 }} \leq \frac{1}{k}
\sum_{i=1}^k e^{a_i}
\label{jensen_numbers}
\end{equation}
Replace $X_i^*$ by $max(X_i^*, 1)$, and combine inequalities
(\ref{k_max_eq}) and (\ref{jensen_numbers}):
$$ \EE e^{\sqrt{ \frac{1}{k} \sum_{i=1}^k (X_i^*)^2 }} \leq 
\EE \frac{1}{k} \sum_{i=1}^k e^{X_i^* + 1} \leq C^{\prime} \frac{n}{k} $$
Another application of Jensen inequality ($\EE \log X \leq \log
\EE X$) yields:
$$ \EE \sqrt{ \frac{1}{k} \sum_{i=1}^k (X_i^*)^2 } \leq \log C^{\prime} \frac{n}{k} $$
which concludes the proof. \hfill $\square$

\medskip \emph{Remark:} If $X_i$ are $\psi_2$ variables, then it can
be simply verified that:
$$ \EE \sqrt{ \frac{1}{k} \sum_{i=1}^k (X_i^*)^2 }\leq c \sqrt
{ \log{ \frac{2n}{k} } } $$ 

\medskip \emph{Proof of Proposition \ref{log_decay}:} Let $\| x \| =
t \| x \|^{\prime}_{\frac{\sqrt{n}}{t}} $, the dual norm of $K_t$. Take
two random bases $\{u_i\}_{i=1}^n$ and $\{v_i\}_{i=1}^n$. The
symmetrized norm $\|| \cdot \||$ is:
$$ \|| x \|| = \EE_{\varepsilon, \varepsilon^{\prime}} \| \sum_{j,k}
\varepsilon_j \varepsilon^{\prime}_k  \langle x, v_j \rangle \langle
v_j, u_k \rangle u_k \| $$
where $\varepsilon, \varepsilon^{\prime}$ are independent and uniformly
distributed in $\{ \pm 1 \}^n$. For $1 \leq i \leq n$ and
$\varepsilon, \varepsilon^{\prime} \in \{ {\pm 1} \} ^n$ define:
$$ \phi_x^i(\varepsilon, \varepsilon^{\prime}) = | \sum_{j,k}
\varepsilon_j \varepsilon^{\prime}_k \langle x, v_j \rangle \langle v_j, u_k
\rangle \langle u_k, e_i \rangle | $$
Then: 
$$ \||x\|| = t \EE_{\varepsilon, \varepsilon^{\prime}}
\big[ \| \phi^1_x( \varepsilon, \varepsilon^{\prime} ),.., \phi^n_x
( \varepsilon,\varepsilon^{\prime} ) \|^{\prime}_{\frac{\sqrt{n}}{t}}
\big] $$ 
By Lemma \ref{convolution_formula}, 
$$ \|| x \|| \leq \sqrt{2} t \EE_{\varepsilon, \varepsilon^{\prime}}
\sqrt { \sum_{i=1}^{\lfloor \frac{n}{t^2} \rfloor +1} \phi^{i*}_x( \varepsilon,
\varepsilon^{\prime})^2 } $$
The following lemma, estimating the $\psi_1$ norm of those variables, will be
proved later.

\begin{lemma}
$$ \| \phi^i_x \|_{\psi_1} < \frac{c_3}{\sqrt n} |x|$$
with probability greater than $1 - e^{-c\sqrt{n}}$ of choosing the
orthogonal bases.
\label{double_random_basis}
\end{lemma}

Lemma \ref{prob_lemma} may now be used (for $k =
\lfloor \frac{n}{t^2} \rfloor +1 $). It shows that:
$$ \|| x \|| \leq \sqrt{2} t (\frac{c_3}{\sqrt n} |x|) \cdot
c_2 (\frac{\sqrt{n}}{t}+1) \log{\frac{2n}{\frac{n}{t^2}} } \leq c
\log t |x| $$

with probability greater than $1 - ne^{-c \sqrt{n}}$ of choosing the
bases. \hfill $\square$

\medskip Before turning to the proof of lemma
\ref{double_random_basis}, prove another lemma, which is believed to
be known to experts:

\begin{lemma}
Let $x = (x_1,..,x_n)$ and $y = (y_1,..,y_n)$ be two random
independent vectors in $S^{n-1}$. Then with probability greater than
$1 - e^{-c \sqrt{n}}$,
$$ \sum_i x_i^2 y_i^2 \leq \frac{c}{n} $$
\label{prob_calculation}
\end{lemma}

\emph{Proof of Lemma \ref{prob_calculation}:} Let $\{ \gamma_i
\}_{i=1}^n$ and $\{ \eta \}_{i=1}^n$ be independent standard
Gaussian variables. Since the measure on the sphere is the radial
projection of the standard Gaussian measure in $\RR^n$, then:
$$ Prob \{ \sum_i x_i^2 y_i^2 > t \} = Prob \{ \frac{1}{\sum_j
\gamma_j^2 \sum_j \eta_j^2} \sum_i \gamma_i^2 \eta_i^2 > t \} $$
To prove the lemma, it is sufficient to bound from below $\sum_j
\gamma_j^2 \sum_j \eta_j^2$ and bound from above $\sum_i
\gamma_i^2 \eta_i^2$. Begin with the second expression. Note that
$\gamma_i^2 \eta_i^2$ is a $\psi_{\frac{1}{2}}$ variable:
$$ (\EE \gamma_i^{2p} \eta_i^{2p} )^{\frac{1}{p}} = (\EE
\gamma_i^{2p}) ^ {\frac{4}{2p}} \leq (c\sqrt{p})^4 = c^4
p^{\frac{1}{\alpha}} $$ 
for $\alpha=\frac{1}{2}$. Therefore, $\sum_i \gamma_i^2 \eta_i^2$ is a
sum of independent copies of a $\psi_{\frac{1}{2}}$ random
variable. By a deviation inequality for sums of i.i.d $\psi_{\alpha}$
random variables (see \cite{S}),
$$ Prob \{ \sum_i \gamma_i^2 \eta_i^2 > cn \} < \text{exp} ( -c^{\prime}
\sqrt{n} ) $$

\medskip The fact that $Prob \{ \sum_j \gamma_j^2 < \frac{n}{2} \}
< e^{-cn}$ follows from Large Deviations technique (e.g. Cram\'{e}r's
Theorem, \cite{V}). To conclude, with probability greater than $1 -
e^{-c \sqrt{n}}$,  
$$ \frac{1}{\sum_j \gamma_j^2 \sum_j \eta_j^2} \sum_i \gamma_i^2
\eta_i^2 < \frac{cn}{\frac{n}{2} \cdot \frac{n}{2}} = \frac{c^{\prime}}{n}
$$

\emph{Proof of Lemma \ref{double_random_basis}:} Let
$\phi_x^i(\varepsilon, \varepsilon^{\prime}) = | \sum_{j,k} 
\varepsilon_j \varepsilon^{\prime}_k \langle x, v_j \rangle \langle
v_j, u_k \rangle \langle u_k, e_i \rangle | $. This random variable
is a particular case of a Rademacher Chaos variable. It is well known
(e.g. see \cite{Tal}), that a $\psi_1$ estimate holds true for such
variables: 
$$\| \phi^i_x \|_{\psi_1} \leq c \| \phi^i_x \|_2 =  c \sqrt{ \sum_j
\langle x,v_j \rangle^2 \sum_k \langle v_j, u_k \rangle^2 \langle u_k,
e_i \rangle^2 } $$ 
\medskip It is sufficient to show that the inequality 
$\sum_k \langle v_j, u_k \rangle^2 \langle u_k, e_i \rangle^2 \leq
\frac{c}{n}$ holds with high probability, since in that case, with the
same probability: 
$$ \sqrt{ \sum_j \langle x,v_j \rangle^2 \sum_k \langle v_j, u_k
\rangle^2 \langle u_k, e_i \rangle^2 } \leq \frac{\sqrt c}{\sqrt n}
\sqrt{ \sum_j \langle x,v_j \rangle^2 } = \frac{\sqrt c}{\sqrt n} |x| $$
The fact that $\sum_k \langle v_j, u_k \rangle^2 \langle u_k, e_i
\rangle^2 \leq \frac{c}{n}$ holds with probability greater than $1 -
e^{-c \sqrt{n}}$ follows directly from Lemma \ref{prob_calculation}:
Take $U \in O(n)$ such that $U(u_k) = e_k$. $U$ is distributed
uniformly over $O(n)$. 
$$ \sum_k \langle v_j, u_k \rangle^2 \langle u_k, e_i
\rangle^2 = \sum_k \langle Uv_j, e_k \rangle^2 \langle e_k, Ue_i
\rangle^2 $$
Since $Uv_j$ and $Ue_i$ are independent and distributed uniformly
over the sphere - the claim is proven, by Lemma
\ref{prob_calculation}. \hfill $\square$

\bigskip \emph{Remark}: Proposition \ref{log_decay} may be adapted
to suit Walsh-type symmetrizations. If $K_t = \sqrt{n} B(l_1^n)
\bigcap t B(l_2^n)$ is symmetrized with respect to Walsh 
vectors $w_1,..,w_{n-1}$, a slightly better conclusion than that in
Proposition \ref{log_decay} is obtained; In this setting, it is 
true that:
$$ \tilde{K_t} \ \subset \ C \sqrt{\log t} B(l_2^n) $$
The differences between the proofs are minor. Lemma
\ref{double_random_basis} becomes much easier as it follows
immediately from Khinchine inequality, even with a $\psi_2$ estimate
rather than $\psi_1$. To take advantage of this improvement, 
use the remark after Lemma \ref{prob_lemma}, to obtain the better
conclusion.

\medskip Re-iteration of this proposition, where each iteration
uses a Walsh basis with respect to the previous, would result in a
rapid decay of the body's diameter. After $\log^* n$ iterations,
a body whose $\frac{ diam(\tilde{K}) }{ M^*( \tilde{K} ) }$ ratio is
bounded by a universal constant is obtained. Note that
this specific choice of symmetrizations decreases the diameter of
all possible convex bodies in $\RR^n$, to be a constant times
their mean width. Of course, once the $\frac{ diam(K) }{ M^*( K )
}$ ratio is bounded, $cn$ random independent Minkowski
symmetrizations suffice for transforming the body into an
approximate Euclidean ball.

\section{Third step: Concentration Techniques}

Take any convex body $K$ in $\RR^n$. According to Corollary
\ref{corollary_2}, from the previous steps (which consist of
$4n$ symmetrizations) a new body is obtained, with  $M^* = 1$ and
with diameter less than $c \log{ \log n}$. As before, the third step
involves symmetrizing with respect to two random orthogonal bases. A
total of $2 n$ symmetrizations will make the body very close to Euclidean. 

\medskip Let $\| \cdot \|$ be the dual norm of the body obtained after
the previous steps. Since $M^*(K) = 1$, then $M(\| \cdot \|) 
\equiv \int_{S^{n-1}} \|x\|d\sigma(x) = 1$, and $b(\| \cdot \|)
\equiv \sup_{x \in S^{n-1}} \|x\| \leq c \log{ \log n}$.
Let $\{u_i\}_{i=1}^n, \{v_i\}_{i=1}^n$ be random orthogonal bases
and denote for $x \in \RR^n$ a set:
$$\FF(x) = \{ \sum_{i,j} \varepsilon_i \varepsilon^{\prime}_j \langle
x, v_i \rangle \langle v_i, u_j \rangle u_j : \varepsilon,
\varepsilon^{\prime} \in \{ \pm 1 \}^n \}$$ 
The symmetrized norm $\|| \cdot \||$ satisfies $\|| x \|| =
\frac{1}{4^n} \sum_{v \in \FF(x)} \| v\|$. This section will prove
that for the new norm: 
$$ \forall x \in \RR^n \ \ \ \ \frac{1}{2} |x| \leq \|| x \|| \leq 2
|x| $$ 
with large probability of choosing $\{u_i\}_{i=1}^n,\{v_i\}_{i=1}^n
\in O(n)$. 
In fact, a somewhat stronger theorem is proved, where instead of
$\frac{1}{2}$ and $2$, better estimates are given.

\bigskip Useful remark: Let $\|| x \|| = \frac{1}{4^n} \sum_{v \in
\FF(x)} \| v\|$ be the norm obtained after symmetrizing with
respect to $\{u_i\}$ and $\{v_i\}$. Take $U \in O(n)$, and let
$\|| \cdot \||_U$ be the norm obtained after symmetrizing with
respect to $\{U u_i \}$ and $\{U v_i\}$. Then $\|| U x \||_U =
\frac{1}{4^n} \sum_{v \in \FF(x)} \| Uv\|$. Therefore, due to the
rotation invariance of the measure $\mu_n$ in $O(n)$, it is possible
to fix an orthonormal system $\{u_i\}$, and prove the following: 

\begin{theorem}
With the above definitions,
$$ \forall x \in S^{n-1} \ \ \  (1 - c\frac{(\log{\log
n})^{\frac{3}{2}}}{\sqrt{\log n}}) \leq \|x\|_U \leq (1 +
c\frac{(\log{\log n})^{\frac{3}{2}}}{\sqrt{\log n}}) $$ 
with probability greater than $1 - e^{-Cn}$ of choosing $U \in
O(n)$, and probability greater than $1 - \frac{1}{n^{10}}$ of
choosing $\{v_i\}$.
\label{6n_final_theorem}
\end{theorem}

The proof shall use three lemmas:

\begin{lemma}
$\forall x, y \in S^{n-1} \ \ \ \frac{1}{4^n} \sum_{v \in \FF(x)}
|\langle v, y \rangle| \leq 2c_1 \frac{\sqrt{\log n}}{\sqrt n} $

for any $\{u_i\}$, with probability of choosing $\{v_i\}$ greater
than $1 - \frac{1}{n^{10}}$.
\label{unbiased_directions}
\end{lemma}

\emph{Proof:} According to Lemma \ref{maximum_of_psi_2}, with
probability greater than $1 - \frac{1}{n^{10}}$, for all $i,j$ the
inequality $|\langle v_i, u_j \rangle| \leq c_1 \frac{\sqrt{\log
n}}{\sqrt n}$ holds. Thus:
$$ \frac{1}{4^n} \sum_{v \in \FF(x)} |\langle v, y \rangle| =
\EE_{\varepsilon, \varepsilon^{\prime}} \big | \sum_{i,j} \varepsilon_i
\varepsilon^{\prime}_j \langle x, v_i \rangle \langle v_i, u_j \rangle
\langle u_j, y \rangle \big | $$
$$ \leq \sqrt{ \EE_{\varepsilon, \varepsilon^{\prime}} \big [ \sum_{i,j}
\varepsilon_i \varepsilon^{\prime}_j \langle x, v_i \rangle \langle
v_i, u_j \rangle \langle u_j, y \rangle \big ]^2 } $$
$$ = \sqrt{ \sum_{i,j} \langle x, v_i \rangle^2 \langle
v_i, u_j \rangle^2 \langle u_j, y \rangle^2 } \leq c_1
\frac{\sqrt{\log n}}{\sqrt n} $$
since $x$ and $y$ are sphere vectors.

\medskip The next lemma is copied from \cite{BLM}, where it is proven.

\begin{lemma}
Assume that $\{ w_\alpha \}_{\alpha \in A} \ \ \subset S^{n-1}$,
and for some $ \delta > 0$,
$$ \sup_{y \in S^{n-1}} \frac{1}{\# A} \sum_{\alpha \in A} |\langle
w_\alpha,y\rangle | \leq \delta $$ 
Let $0 < \lambda < 1$ and let $k \leq n$ be an integer. Then the
set $A$ can be partitioned into families $\mathcal{F}_\beta = \{
\beta_i \}_{i=1}^k \subset A, \ \beta \in B$, so that
$\# (\cup_{\beta \in B} \mathcal{F}_\beta) > (1-\lambda)\#A - k$
and so that for every $\beta \in B$ there is an orthonormal set of
vectors $\{ v_{\beta_i} \}_{i=1}^k$ satisfying
$$ |v_{\beta_i} - w_{\beta_i}| \leq \frac{\delta 4^k}{\lambda} $$
\label{combinatorial_lemma}
\end{lemma}

Concentration on the orthogonal group shall be used in the proof of
Theorem \ref{6n_final_theorem}, due to \cite{GM} (see \cite{MS}, page 29): 

\begin{lemma}
Let $\| \cdot \|$ be a norm on $\RR^n$ such that $\|x\| \leq b|x|
\ \forall x \in S^{n-1}$. Let $k \leq n$ be a positive integer,
and $\{x_i\}_{i=1}^k$ be orthonormal vectors. Denote $M =
\int_{S^{n-1}} \|x\| d\sigma_n(x)$. Then:
$$ \mu_n \{ U \in O(n) \ ; \ | \frac{1}{k} \sum_{i=1}^k \|Ux_i\| - M| \geq
\varepsilon \} \leq \text{exp} (-c_4 \frac{\varepsilon^2 nk}{b^2})$$
\label{ortho_concentration}
\end{lemma}

\bigskip \emph{Proof of Theorem \ref{6n_final_theorem}}:
Fix $x \in S^{n-1}$. Let $\varepsilon = c_5 \frac{(\log{ \log
n})^{\frac{3}{2}}}{\sqrt{\log n}}$, $\lambda = \frac{\varepsilon}{b}$ and $k =
\frac{\log n}{10}$. According to Lemma \ref{unbiased_directions}, the
collection of vectors $\FF(x)$ satisfies the requirement of Lemma
\ref{combinatorial_lemma} for $\delta = 2c_1 \frac{\sqrt{ \log
n}}{\sqrt n}$, with large probability of choosing $\{v_i\}$ (and of
course independently of $U$). As a result, $\FF(x)$ can be
decomposed into disjoint almost orthogonal families $\{ \FF_\beta
\}_{\beta \in B}$, which cover all but a $\lambda$ fraction of
$\FF(x)$. \\
From Lemma \ref{combinatorial_lemma}, for each family $\FF =
\{x_1,..,x_k\} \subset \FF(x)$, there exist orthonormal vectors
$\{t_1,..,t_k\}$ such that $|t_i - x_i| \leq \frac{\delta
4^k}{\lambda}$. Since $\{t_i\}_{i=1}^k$ are orthonormal, then by
Lemma \ref{ortho_concentration}:
$$ \mu_n \{ U \in O(n) \ ; \ | \frac{1}{k} \sum_{i=1}^k \|Ut_i\| - 1| \geq
\varepsilon \} \leq \text{exp} (-c_4 \frac{nk \varepsilon^2}{b^2}) $$
where $b = \sup_{x \in S^{n-1}} \|x\|$. Since $\|Ut_i - Ux_i\| \leq b
\frac{\delta 4^k}{\lambda}$, then:
\begin{equation}
|\frac{1}{k} \sum_{i=1}^k \|Ux_i\| - 1 | \leq \varepsilon + b
\frac{\delta 4^k}{\lambda} \label{equation_1}
\end{equation}
with probability (of choosing $U \in O(n)$) of at least $1 -
\text{exp} (-c_4 \frac{nk \varepsilon^2}{b^2})$.

\medskip This holds for a single family $\FF$. The number of families is
less than $4^n$, so inequality (\ref{equation_1}) holds for all
families $\{ \FF_\beta \}_{\beta \in B}$ together, with
probability greater than $1 - 4^n \text{exp} (-c_4 \frac{nk
\varepsilon^2}{b^2}) = 1 - \text{exp} (-c_4 n (\frac{k
\varepsilon^2}{b^2}-\log 4) ) $. \\
There still remains a $\lambda$ fraction of the collection
$\FF(x)$, not covered by the disjoint families $\{
\FF_\beta \}_{\beta \in B}$. Their contribution to the relevant
expression, which is $|\frac{1}{4^n} \sum_{v \in \FF(x)} \| U v \| -
1|$, can be bounded by $\lambda b$. Hence:
$$ | \ \|| Ux \||_U - 1 \ | = |\frac{1}{4^n} \sum_{v \in
\FF(x)} \| Uv \| - 1| \leq $$
$$ \leq \frac{k}{4^n} \sum_{\beta \in B}
|\frac{1}{k} \sum_{v \in \FF_\beta} (\|Uv\|-1)| + \frac{\lambda
4^n}{4^n} b \leq (1-\lambda)(\varepsilon + b \frac{\delta
4^k}{\lambda}) + \lambda b $$

\bigskip In summary: choose $\{v_i\}$ by random. With probability of
at least $1 - \frac{1}{n^{10}}$, the following holds: the set of $U
\in O(n)$ for which 
\begin{equation}
| \ \|| Ux \||_U - 1 \ | \leq \varepsilon + b \frac{\delta
4^k}{\lambda} + \lambda b \label{eq_2}
\end{equation}
has measure of at least $1 - \text{exp} (-c_4 n (\frac{k \varepsilon^2}{b^2}-
\log 4))$. From substituting the values of the variables $k$,
$\varepsilon$, $\lambda$, it follows that $\lambda b \leq
\varepsilon$, and also $b \frac{\delta 4^k}{\lambda} < \varepsilon$,
for $n > c_6$. Therefore, the quantity discussed in (\ref{eq_2})
is less than $3 \varepsilon$, for $n > c_6$. \\
The inequality $ | \ \|| Ux \||_U - 1 \ | \leq 3 \varepsilon$
holds with probability (with respect to $U$) of at least $1 -
\text{exp} (-c_4 n (\frac{k \varepsilon^2}{b^2}- \log 4))$. With a suitable
universal constant $c_5$ this probability would be greater than $1
- \text{exp} (-10 n \log {\log n})$.

\bigskip This analysis considered a fixed $x \in S^{n-1}$. Now,
take an $\varepsilon$-net on the sphere denoted by $\mathcal{N}$.
There exists such a net with $\# \mathcal{N} \leq
(\frac{4}{\varepsilon})^n$. For each $x \in \mathcal{N}$, $| \ \|| Ux
\||_U - 1 \ | \leq 3 \varepsilon$ with probability greater than $1 - \text{exp}
(-10 n \log {\log n})$. Since $(\frac{4}{\varepsilon})^n \leq
\text{exp} (\log {\log n})$ for $n > c_6$, then $| \ \|| Ux \||_U - 1
\ | \leq 3 \varepsilon$ holds for all $x \in \mathcal{N}$, with more
than exponentially close to $1$ probability.

\medskip For a general $x \in S^{n-1}$, write $x =
\sum_{i=0}^{\infty} \theta_i x_i$, where $Ux_i \in \mathcal{N}$,
and $\theta_0 = 1, 0 \leq \theta_i \leq \varepsilon^i$. Then $\|| x
\||_U \leq \sum_{i=0}^{\infty} (1+3 \varepsilon) \varepsilon^i =
\frac{1+3\varepsilon}{1 - \varepsilon} \leq 1 + 5\varepsilon$. Finally,
$\||x\||_U \geq \||x_0\||_U - \sum_{i=1}^{\infty} |\theta_i| \cdot
\||x_i\||_U \geq 1 - 5 \varepsilon$.

\bigskip Hence, with slightly better than exponentially close to $1$
probability, the new norm $\|| \cdot \||$ satisfies
$$ \forall x \in \RR^n \ \ \ (1-\varepsilon)|x| \leq \||x\|| \leq
(1+\varepsilon)|x| $$
where $\varepsilon < c \frac{(\log {\log n})^{\frac{3}{2}}}{\sqrt{\log
n}}$, and the theorem is proven, for $n > c_6$. \hfill $\square$

\medskip \emph{Remark:} Using a Walsh-type symmetrization in the
second step, the theorem can be proven with $\varepsilon < c
\frac{\log {\log n}}{\sqrt{\log n}}$, an improvement of a mere
$\sqrt{ \log{\log n} }$ factor. 

\bigskip I would like to express my sincere thanks and appreciation
to my supervisor, Prof. Vitali Milman, for our inspiring
discussions and for his ongoing support throughout the research of
this topic.

\end{document}